\documentclass[12pt]{article}

\usepackage{amssymb}

\newtheorem{theorem}{Theorem}

\newtheorem{lemma}[theorem]{Lemma}

\newenvironment{remark}[1][Remark]
           {\medbreak\noindent \textbf{#1 \enspace}}
           {\par \medbreak}

\newenvironment{proof}[1][Proof]
           {\medbreak\noindent \emph{#1: \enspace}}
           {\hfill $\dashv$  \par \medbreak}

\newcommand{\bproof}[1][Proof]
           {\medbreak\noindent \emph{#1: \enspace}}
\newcommand{\eproof}[1]
           {\hfill$\dashv$~#1  \par \medbreak}

\begin{document}
\baselineskip=18pt

\thispagestyle{empty}

\title{New $\Sigma^1_3$ Facts}

\author{Sy D. Friedman\thanks{%
  Research supported by NSF Contract \#9625997-DMS}\\[1ex]
  {\normalsize M.I.T.}
}
\date{}

\maketitle

\thispagestyle{empty}

\noindent
{\bf Abstract} We use ``iterated square sequences'' to show: There is an $L$-definable 
partition $n:L-Singulars \to \omega$ such that if $M$ is an inner model 
without $0^\#$:
\noindent
(a) For some $n, M \models \{\alpha|n(\alpha)\leq n\}$ is stationary.
\noindent
(b) For each $n$ there is a generic extension of $M$ in which $0^\#$ does
not exist and $\{\alpha|n(\alpha)\leq n\}$ is non-stationary.
\noindent
This result is then applied to show that if $M$ is an inner model
without $0^\#$ then some $\Sigma^1_3$ sentence not true in $M$ can be 
forced over $M$. 
\vspace {2ex}

\noindent
Assume that $0^\#$ exists and that $M$ is an inner model of ZFC,
$0^\# \notin M$.  Then of course $M$ is not $\Sigma^1_3$-correct:
the true $\Sigma^1_3$ sentence ``$0^\#$ exists'' is false in~$M$.
In this article we use a result about $L$-definable partitions
(which may be of independent interest) to show that in fact this
effect can be achieved by forcing over~$M$.  We work in
Morse-Kelly class theory.

\begin{theorem}
  \label{th:1}
  Assume that $0^\#$ exists.  There exists an $\omega$-sequence of true
  $\Sigma^1_3$ sentences $\langle \varphi_n \mid n \in \omega \rangle$
  such that if $M$ is an inner model, $0^\# \notin M$:

  \renewcommand{\theenumi}{\alph{enumi}}%
  \renewcommand{\labelenumi}{(\theenumi)}%
  \begin{enumerate}
  \item 
    $\varphi_n$ is false in $M$ for some $n$.

  \item 
    For each $n$, some generic extension of $M$ satisfies
    $\varphi_n$.  
  \end{enumerate}

  Moreover if $M = L [R]$, $R$ a real then these
  generic extensions can be taken as inner models of $L [R, 0^\#]$.
  
\end{theorem}

The above result is based on the next result, concerning
$L$-definable partitions.

\begin{theorem}
  \label{th:2}
  There exists an $L$-definable function $n: L \hbox{-Singulars}
  \to \omega$ such that if $M$ is an inner model, $0^\# \notin M$:

  \renewcommand{\theenumi}{\alph{enumi}}%
  \renewcommand{\labelenumi}{(\theenumi)}%
  \begin{enumerate}
  \item 
    For some $n$, $M \models \{ \alpha \mid n (\alpha) \le n \}$
    is stationary.

  \item 
    For each $n$ there is a generic extension of $M$ in which
    $0^\#$ does not exist and $\{ \alpha \mid n (\alpha) \le n \}$
    is non-stationary.
  \end{enumerate}
\end{theorem}

\begin{remark}
  ``Stationary in $M$'' means: intersects every $M$-definable
  (with parameters) $CUB$.
\end{remark}

\begin{proof}
  We define $n (\alpha)$.  Let $\langle C_\alpha \mid \alpha \; L
  \hbox{-singular} \rangle$ be an $L$-definable
  $\square$-sequence: $C_\alpha$ is $CUB$ in $\alpha$, $otC_\alpha =$
  ordertype $C_\alpha < \alpha$ and $\bar{\alpha} \in \lim
  C_\alpha \to C_{\bar{\alpha}} = C_\alpha \cap
  \bar{\alpha}$.  If $otC_\alpha$ is $L$-regular then $n
  (\alpha) = 0$.  Otherwise $n (\alpha) = n (otC_\alpha) + 1$.

  (a) is clear, as otherwise there is a $CUB$ $C \subseteq
  L$-regulars amenable to $M$, contradicting that Covering
  Theorem and the hypothesis that $0^\#$ does not belong to~$M$.

  Now we prove (b).  Fix $n \in \omega$.  In $M$ let $P$ consist of
  closed, bounded $p \subseteq$ ORD such that $\alpha \in p \to
  \alpha$ $L$-regular or $n (\alpha) \ge n + 1$, ordered by $p
  \le q$ iff $p$ end extends~$q$.
  
  We claim that $P$ is $\infty$-distributive in~$M$.  Suppose
  that $p \in P$ and $\langle D_\alpha \mid \alpha < \kappa
  \rangle$ is a definable sequence of open dense subclasses of
  $P$, $\kappa$ regular.  We wish to find $q \le p$, $q \in
  D_\alpha$ for all $ \alpha < \kappa$.  Let $C = \{ \beta \mid
  \beta \hbox{ a strong limit cardinal, for all } \alpha <
  \kappa: r \in V_\beta \to \exists s \le r (s \in V_\beta, \; s
  \in D_\alpha) \}$, a $CUB$ class of ordinals.  It suffices to
  show that $C \cap \{ \beta \mid n (\beta) \ge n + 1 \}$ has a
  closed subset of ordertype $\kappa+1$, for then $p$ can be
  successively extended $\kappa$ times meeting the $D_\alpha$'s,
  to conditions with maximum in $\{ \beta \mid n (\beta) \ge n +
  1 \}$; the final condition (at stage $\kappa$) extends $p$ and
  meets each $D_\alpha$.

  \begin{lemma}
    \label{lem:3}
    Suppose $m \ge n$, $\alpha$ is regular and $C$ is a closed
    set of ordinals greater than $\alpha^{+ m}$ of ordertype
    $\alpha^{+ m} + 1$ (where $\alpha^{+0} = \alpha$, $\alpha^{+
      (k + 1)} = (\alpha^{+k})^+$).  Then $C \cap \{ \beta \mid
    n(\beta) \ge n \}$ has a closed subset of ordertype
    $\alpha^{+ (m - n)} + 1$.
  \end{lemma}

  \bproof[Proof of Lemma~\ref{lem:3}] 
    By induction on $n$.  Suppose $n = 0$.  Let $\beta = \max C$.
    Then $\beta$ is singular and hence singular in $L$.  So
    $C_\beta$ is defined and $\lim (C_\beta \cap C)$ is a closed
    set of ordertype $\alpha^{+m} +1$ consisting of
    $L$-singulars.  So $\lim (C_\beta \cap C) \subseteq C \cap \{
    \gamma \mid n (\gamma) \ge 0 \}$ satisfies the lemma.

    Suppose the lemma holds for $n$ and let $m \ge n$, $C$ a
    closed set of ordertype $\alpha^{+ (m+1)} +1$ consisting of
    ordinals greater than $\alpha^{+ (m+1)}$.  Let $\beta = \max
    C$.  Then $ C_\beta$ is defined and $D = \lim (C_\beta \cap
    C)$ is a closed set of ordertype $\alpha^{+ (m+1)} +1$.  Let
    $\bar{\beta} = (\alpha^{+m} + \alpha^{+m} + 1)$st element of
    $D$.  Then $\bar{D} = \{ otC_\gamma \mid \gamma \in D$,
    $(\alpha^{+m} + 1)\hbox{st element of } D \le \gamma \le
    \bar{\beta} \}$ is a closed set of ordertype $\alpha^{+m} +
    1$ consisting of ordinals greater than $\alpha^{+ m}$.  By
    induction there is a closed $\bar{D}_0 \subseteq \bar{D} \cap
    \{ \gamma \mid n (\gamma) \ge n \}$ of ordertype $\alpha^{+
      (m-n)} + 1$.  But then $D_0 = \{ \gamma \in D \mid
    otC_\gamma \in \bar{D}_0 \}$ is a closed subset of $C \cap \{
    \gamma \mid n (\gamma) \ge n + 1 \}$ of ordertype $\alpha^{+
      (m-n)} + 1$.  As $\alpha^{+ (m-n)} = \alpha^{+ ((m+1) -
      (n+1))}$ we are done.
  \eproof{(Lemma~\ref{lem:3})}

  By the lemma, $C \cap \{ \beta \mid n (\beta) \ge n \}$ has
  arbitrary long closed subsets for any $n$, for any $CUB$ $C
  \subseteq$ ORD.  It follows that $P$ is $\infty$-distributive.
  Now to prove (b), we apply the forcing $P$ to $M$, producing
  $C$ witnessing the nonstationarity of $\{ \alpha \mid n
  (\alpha) \le n \}$, and then follow this with the forcing to
  code $\langle M, C \rangle$ by a real, making $C$ definable.
  Of course this will not produce $0^\#$ as every successor to a
  strong limit cardinal is preserved in the coding.
\end{proof}

We also note that in Theorem~\ref{th:2} the generic extension can
be formed in $L [R, 0^\#]$ in the case $M = L [R]$, $R$ a real,
using the fact that in $L [R, 0^\#]$, generics can be constructed
for $P$ (an ``Amenable'' forcing) and for Jensen coding
(see~\cite[Friedman]{99}).

\begin{proof}[Proof of Theorem~\ref{th:1}]
  We use David's trick (see~\cite[Friedman]{98}).  Let $\varphi_n$
  be the $\Sigma^1_3$ sentence: $\exists R \forall \alpha (
  L_\alpha [R] \models Z F^{-} \to L_\alpha [R] \models \beta
  \hbox{ a limit cardinal } \to \beta \, L\hbox{-regular or } n
  (\beta) \ge n)$.  By Theorem~\ref{th:2}(b) and cardinal
  collapsing (to guarantee that limit cardinals $\beta$ are either
  $L$-regular or satisfy $n (\beta) \ge n$), $M$ has a generic extension
  $L[R] \models \beta \hbox { a limit cardinal} \to 
  \beta \, L\hbox{-regular or }
  n(\beta)\ge n$ (inside $L [S, 0^\#]$ if $M =
  L [S]$, $S$ a real).  By David's trick we can in fact obtain
  $\varphi_n$ in $L [R]$.
\end{proof}

\begin{remark}[Question]
  Can the generic extensions in Theorem~\ref{th:1}(b) be taken to
  have the same cofinalities as $M$, in case $M$ satisfies $G C H$?
\end{remark}

\end{document}